\documentclass[12pt]{article}
\usepackage[british]{babel}
\usepackage{mathrsfs}
\usepackage{amsmath}
\usepackage{amssymb}
\usepackage{amsthm}
\usepackage{enumerate}
\usepackage{epic}
\usepackage{color}

\numberwithin{equation}{section}
\makeatletter
\renewcommand{\subsection}{\@startsection
{subsection}{2}{0mm}{\baselineskip}{-0.25cm}
{\normalfont\normalsize\bf}}
\makeatother

\newtheorem{theorem}{Theorem}[section]
\newtheorem{proposition}[theorem]{Proposition}
\newtheorem{lemma}[theorem]{Lemma}
\newtheorem{corollary}[theorem]{Corollary}

\newtheorem{result}[theorem]{Result}

\def\F{\mathbf F}

\def\cF{\mathcal F}

\def\cX{\mathcal X}

\def\xfq2{{\cX(\mathbb{F}_{q^2})}}
\def\F+xfq2{{\cF^+(\mathbb{F}_{q^2})}}
\def\xpfq2{{\cX^+(\mathbb{F}_{q^2})}}
\def\xmfq2{{\cX^-(\mathbb{F}_{q^2})}}
\title{Circles and Paths in $2$-Colored Best Match Graphs }
\author{Annachiara Korchmaros\footnote{Annachiara Korchmaros: Department of Mathematics and Systems Analysis, Aalto University, annachiara.korchmaros@aalto.fi}}
\date{}

\begin{document}
\date{}
\maketitle
\begin{abstract}
Recent investigations in computational biology focus on a family of $2$-colored digraphs, called $2$-colored best match graphs, which naturally arise from rooted phylogenetic trees. Actually the defining properties of such graphs are unexpectedly unusual in graph theory, and they were established only recently after the discovery of their links to evolutionary relatedness via phylogenetic trees. In this paper several results are obtained on $2$-colored best match graphs which well fit in the mainstream of graph theory.
\end{abstract}

\section{Introduction}
Best match graphs have been an important tool in current studies of computational biology~\cite{geiss,geiss1,geiss3}. The formal definition of such graphs comes from evolutionary relatedness via phylogenetic trees. Let $T$ be a rooted phylogenetic tree with leaf set $L$ and a surjective color-map $\sigma:L\rightarrow S$ for a non-empty color set $S$. Then $y\in L$ is a best match of $x\in L$, in symbols $x\rightarrow y$, if ${\rm{lca}}(x,y)\prec {\rm{lca}}(x,y')$ for all $y'\in L$ with $\sigma(y)=\sigma(y')$. Here ${\rm{lca}}$ stands for the last common ancestor, and the partial ordering $p\prec q$ occurs if $q$ is located above $p$ along the path connecting $p$ to the root of $T$. The associated colored best match graph, shortly cBMG, is the directed graph on the vertex set $L$ where arcs are the ordered pairs $xy$ with $x\rightarrow y$ and  $x\neq y$. In particular, any cBMG is a colored digraph with color map $\sigma$. If a colored digraph $G(T,\sigma)$ is isomorphic to the ${\rm{cBMG}}$ associated to the rooted phylogenetic tree $T$, then $T$ is said to explain the vertex-colored graph $(G,\sigma)$. In any ${\rm{cBMG}}$ there is a natural equivalence relation $\dot{\sim}$ where $x{\dot{\sim}}y$ if $x$ and $y$ have the same out-neighbours and in-neighbours. In particular, if $N(x)$ denotes the set of all out-neighbours of a vertex $x$, then $x\dot{\sim}y$ implies $N(x)=N(y)$ but the converse is not always true. In fact, if $N^-(x)$ denotes the set of all in-neighbours of a vertex $x$ then  $N^-(x)\neq N^-(y)$ may occur even though $N(x)=N(y)$ holds.

In the case where $|S|=2$, ${\rm{cBMG}}$ is called a 2-cBMG which is a bipartite digraph. The study of 2-cBMGs is of special interest since all cBMGs have some induced 2-cBMG subgraphs; see \cite[Theorem 9]{geiss}. The fundamental properties of connected 2-cBMGs are found in \cite{geiss}. In particular, \cite[Theorem 3]{geiss}  
establishes that for any two vertices $u$ and $v$ belonging to different equivalence classes, the following properties hold. 
\begin{itemize}
\item[N(1)] $u\cap N(v)=v\cap N(u)=\emptyset$ implies $N(u)\cap N(N(v))=N(v)\cap N(N(u))=\emptyset.$
\item[N(2)] {\mbox{$ N(N(N(u))) \subseteq N(u)$.}}
\item[N(3)]  $u\cap N(N(v))=v\cap N(N(u))=\emptyset$ together with $N(u)\cap N(v)\neq \emptyset$ implies $N^-(u)=N^-(v)$ and one of the inclusions $N^-(u)\subseteq N^-(v)$, $N^-(v)\subseteq N^-(u)$.
\end{itemize}
Actually, the above three properties are sufficient to characterize connected 2-cBMGs. In fact, for a connected $2$-colored digraph, \cite[Theorem 4]{geiss} states that there exists a rooted phylogenetic tree $T$ explaining $(G,\sigma)$ if and only if $(G,\sigma)$ satisfies $N(1),N(2)$, and $N(3)$ for any two vertices $u$ and $v$ from two different equivalence classes.

It should be noticed however that $N(1)$ and $N(3)$ are properties which had not been considered at all in the literature on graph theory until the discovery of their links to evolutionary relatedness via phylogenetic trees whereas $N(2)$, under the name of bitransitive property, was marginally studied in a recent manuscript \cite{das} about bitournaments; see Section \ref{bitour}. 

This gives a motivation to consider graph-theoretic properties of ``classical type'' in digraphs satisfying at least one of the above three conditions.
We are mostly concerned with $N(2)$, and our contributions are stated and proven in Section \ref{secpc}. It turns out that our results on circles and paths fit in well with classical works in graph theory dating back to 1980's. In Section \ref{secdis} we point out that the quotient graph of a digraph with property $N(2)$ has a long path. If we assume that
$N(1)$ and $N(3)$ also hold then our results can be refined, see Section \ref{n1n2} and \ref{n2n3}. However it remains unclear how $N(3)$ can affect the structure of the underlying 2-color digraph.

\section{Notation and Terminology}
In this paper $\Gamma$ stands for a digraph without loops, parallel edges, and vertices which have no out-neighbours. With the usual notation, $V$ is its vertex-set, $E$ is its edge-set, and $\Gamma=\Gamma(V,E)$.

For a vertex $u$ of $\Gamma$, $v$ is an \emph{out-neighbor} (respectively \emph{in-neighbor}) of $u$ if $uv$ (respectively $vu$) is an edge of $\Gamma$. The set of all out-neighbors (respectively in-neighbors) of $u$ is denoted by $N(u)$ and $N^-(u)$ respectively. Therefore, $u$ and $v$ have no common out-neighbor (respectively in-neighbor) if and only if $N(u)\cap N(v)=\emptyset$ (respectively  $N^-(u)\cap N^-(v)=\emptyset$). Furthermore, for any two vertices $u$ and $v$ of $\Gamma$, we say that ``$u$ is \emph{out-dominated} by $v$'' if $N(u)\subseteq N(v)$.

For an ordered pair of vertices $(u,v)$ of $\Gamma$, $v$ is \emph{strongly connected} to $u$, if there exists a directed walk $u=u_0\rightarrow u_1 \rightarrow u_2 \rightarrow \cdots\rightarrow u_m \rightarrow u_{m+1}=v $ such that $u_iu_{i+1}$ is an edge of $\Gamma$ for every $0\le i \le m$. The \emph{length} of the directed walk is the number of edges in it, i.e. $m+1$. A \emph{directed trail} is a directed walk in which all edges are distinct. A \emph{directed path} is a directed trail in which all vertices are distinct. A \emph{directed circuit} is a non-empty directed trail in which the first and last vertices are repeated. A \emph{directed cycle} is a directed trail in which all vertices but the first and last are distinct. Two vertices $u$ and $v$ in $\Gamma$ are \emph{independent} if $uv\not \in E$ and $vu\not \in E$, i.e. $v\not \in N(u)$ and $u\not \in N(v)$. It may be noticed that any two equivalent vertices $u$ and $v$ of $\Gamma$ are independent, otherwise $u\in N(v)=N(u)$ would imply that $uu\in E$ which is impossible as $\Gamma$ has no loops.

A \emph{bipartite digraph} $\Gamma(U,V,E)$  is a digraph whose vertices can be divided into two disjoint sets $U$ and $V$ such that every edge connects a vertex in $U$ to one in $V$.
The two sets $U$ and $V$ may be thought of as a coloring of the graph with two colors where a \emph{coloring} is a labeling of the vertices with the two colors such that no two vertices sharing the same edge have the same color. A bipartite graph is \emph{balanced} if $|U|=|V|$.

A digraph is \emph{oriented} if $uv\in E$ implies $vu\not\in E$. An oriented bipartite digraph $\Gamma(U,V,E)$ is \emph{bitransitive} if for all vertices $x_1,x_2,y_1,y_2\in U\cup V$ with $x_1y_1,y_1x_2,x_2y_2 \in E$ we have $x_1y_2 \in E$. An oriented bipartite digraph $\Gamma(U,V,E)$ is a \emph{bitournament} if for all $u\in Y$ and $v\in V$, either $uv\in E$ or $vu\in E$.

For a pair $(\mathcal{A},\mathcal{O})$ where $\mathcal{A}$ is a finite set of non-negative even integers and $\mathcal{O}$ is a set a positive odd integers, the associated \emph{odd-even} oriented digraph ${\vec{\mathcal{G}}}_{\mathcal{A}}(\mathcal{O})$ on vertex-set $\mathcal{A}$ has edge-set $E$ with $ab\in E$ when both $\frac{1}{2}(a+b)$ and $\frac{1}{2}(b-a)$ belong to $\mathcal{O}$. The odd-even graph ${\vec{\mathcal{G}}}_{\mathcal{A}}(\mathcal{O})$ is an oriented bipartite graph $\Gamma(U,V,E)$ with $U=\{a|a\equiv 0 \pmod{4},a\in \mathcal{A}\}$ and
$V=\{a| a\equiv 2 \pmod{4}|,a\in \mathcal{A}\}$.

\section{Cycles and paths in bipartite digraphs}
\label{seccp}
Several papers gave sufficient conditions for bipartite digraphs, in terms of the number of edges, to have cycles and paths with specified properties.
These conditions can be viewed as digraph versions or variants of similar conditions on bipartite graphs which were widely studied since the 1980's; see \cite{am,am1,aj,ch,das,BJ,MMMM,MM,NV,wang,zang1,zang2}. We recall those which are related to the present investigation. The references are \cite{am,am1,ch,MMMM,zang1}.
\begin{result}
\label{resCMM} Let $\Gamma(U,V,E)$ be a bipartite digraph with $|U|=a,|V|=b$, $a\le b$, and $k=\min\{N(x)+N^-(x)|x\in U\cup V\}$.
\begin{itemize}
\item[(i)] If $|E|\ge 2ab-b+1$ then $\Gamma(U,V,E)$ has a cycle of length $2a$.
\item[(ii)] If $|E|\ge 2ab-(k+1)(a-k)+1$ then $\Gamma(U,V,E)$ has a cycle of length $2a$.
\item[(iii)] If $|E|\ge 2ab-k(a-k)+1$ then, for any two vertices $x$ and $y$ which are not in the same partite set, there is a path from $x$ to $y$ of length $2a-1.$
\item[(iv)] If $|E|\ge 2ab-a+2$, then for $x,y\in U\cup V$, any set of $a-1$ vertices is contained in a path of length at least $2(a-1)$ from $x$ to $y$ while for $u\in U, v\in V$, there are paths from $u$ to $v$ and from $v$ to $u$ of every odd length $m$ with $3\le m \le 2a-1$.
\item[(v)] If $a\le 2k-1$ then $\Gamma(U,V,E)$ has a cycle of length $2a$, unless either $b>a=2k-1$ and $\Gamma(U,V,E)\cong\Gamma_1(a,b)$ or $k = 2$ and $\Gamma(U,V,E)\cong \Gamma_2(3,b)$.
\end{itemize}
\end{result}
\begin{result}
\label{zanghA} Let $\Gamma(U,V,E)$ be a bipartite oriented digraph whose vertex in-degree is at least $h\geq 0$ and out-degree is at least $k\geq 0$ for all vertices. Then  $\Gamma(U,V,E)$ contains either a directed cycle of length at least $2(k+h)$ or a directed path of length at least $2(k+h)+3$. 
\end{result}
\begin{result}
\label{wangth}  Let $\Gamma(U,V,E)$ be a balanced directed bipartite graph with $|U|=|V|=n\ge 2$. Suppose that $N(u)+N^-(u)+N(v)+N^+(v)>3n + 1$ for all $u\in U,v\in V$. Then $\Gamma=(U,V,E)$ contains two vertex-disjoint directed cycles of lengths $2n_1$ and $2n_2$, respectively, for any positive integer partition $n = n_l + n_2$.
\end{result}
\begin{result}
\label{dastheorem3.4}
Every acyclic oriented bipartite graph is isomorphic to some odd-even graph.
\end{result}

\section{Bitransitive Bitournaments}\label{bitour}
A nice example of a bitransitive bitournament arises from arithmetic, see \cite{das}. For a nonempty subset $S$ of natural numbers, let $\Gamma_S(S,E)$ be the digraph with vertex set $S$ such that $uv\in E$  if $u<v$ and $u$ and $v$ have opposite parity. If $U$ consists of all even numbers in $S$ while $V$ consists of all odd numbers in $S$ then $\Gamma_S(S,E)$ is bitransitive bitournament. The importance of this example is due to the following characterization; see \cite[Theorem 2.5]{das}.
\begin{result}
\label{dastheorem2.5} Let $\Gamma(U,V,E)$ be a bitournament. Then the following properties are equivalent.
\begin{itemize}
\item[(i)] $\Gamma(U,V,E)$ is bitransitive.
\item[(ii)] $\Gamma(U,V,E)$ has no directed cycle.
\item[(iii)] $\Gamma((U,V,E)\cong \Gamma_S(S,E)$ for some nonempty subset $S$ of natural numbers.
\end{itemize}
\end{result}

\begin{result}
\label{dastheorem2.9} A bitournament $\Gamma(U,V,E)$ with $|U| = |V| = 2m$  and $|N(u)|=|N^-(u)|$ for all $u\in U\cup V$ is not bitransitive.
\end{result}

\section{Quotient graph}
The equivalence relation $\dot{\sim}$ gives rise to the \emph{quotient graph}  $\bar{\Gamma}(\bar{V},\bar{E})$ whose vertices are the equivalence classes and edges are defined as follows. Let $\bar{u}$ and $\bar{v}$ be two vertices of $\bar{\Gamma}$, then the ordered pair $\bar{u}\bar{v}$ is in $\bar{E}$ whenever $uv\in E$ for every $u\in \bar{u}$ and $v\in\bar{v}$.
In other words, $\bar{v}\in N(\bar{u})$ if and only if $v\in N(u)$ for every $u\in \bar{u}$ and $v\in\bar{v}$.
Several properties of $\Gamma$ are inherited by $\bar{\Gamma}$ such as connectivity and vertex colorability.

\begin{lemma}
\label{lem01} $\Gamma$ is connected if and only if $\bar{\Gamma}$ is connected.
\end{lemma}
\begin{proof} We limit ourselves to show that connectivity of $\bar{\Gamma}$ implies the connectivity of $\Gamma$, the converse being trivially true. For $\bar{u},\bar{v}\in \bar{V}$, let
$\bar{u}=\bar{x}_0\rightarrow\bar{x}_1\rightarrow\cdots\rightarrow \bar{x}_m \rightarrow\bar{x}_{m+1}=\bar{y}$ be a walk in the undirected graph $\bar{\Gamma}^*$ associated to $\bar{\Gamma}$.
For every $0\le i\le m+1$, take just one vertex from $\bar{x}_i$, say $x_i$. Then $x_ix_{i+1}$ with $i=0,1,\ldots,m$ is an edge in the undirected graph $\Gamma^*$ associated to $\Gamma$.
Let $u=x_0$ and $v=x_{m+1}$. Then the sequence ${u}={x}_0\rightarrow{x}_1\rightarrow\cdots\rightarrow {x}_m \rightarrow{x}_{m+1}={y}$ is a walk in $\Gamma^*$ which proves the claim.
\end{proof}

The \emph{chromatic number} of a directed graph is the minimum number of colors necessary to color the vertices such that there is no edge incident vertices with the same color.
\begin{lemma}
\label{lem02} $\Gamma$ and $\bar{\Gamma}$ have the same chromatic number.
\end{lemma}
\begin{proof} Let consider an minimal coloring of $\Gamma$ and choose one vertex from each equivalent class of $\Gamma$. Suppose there is a equivalent class $\alpha$ containing a vertex $v$ whose color is different from the one of the vertex $u$ we have chosen in $\alpha$. Now we show that replacing the color of $v$ with the color of $u$ still defines a vertex coloring in $\Gamma$. By way of a contradiction, there is a vertex $w$ with same color of $v$ such that either $vw \in E$ or $wv \in E$. Here $w\neq u$ since $u$ and $v$ are in the same class hence they are independent. Therefore $u$ and $w$ are two different vertices which have the same color since $u$ and $v$ are in the same equivalent class, contradicting our hypothesis. This shows that $\Gamma$ has a minimal coloring such that vertices in the same class must have the same color. This coloring naturally induces a coloring on $\bar{\Gamma}$ with same number of colors.

On the other hand, given a minimal coloring of $\bar{\Gamma}$, one can define a coloring on $\Gamma$ by assigning to $u\in V$ the color of the equivalent class where it belongs to in $\bar{\Gamma}$. Now we need to prove that if $uv\in E$ then $u$ and $v$ have different colors. Let $uv\in E$, two cases are distinguished according as $u$ and $v$ are equivalent or not. We can easily rule out the former case because equivalent vertices are independent. In the latter case, $\bar{u}$ and $\bar{v}$ are two distinct vertices in $\bar{\Gamma}$ and hence those vertices have different colors, then $u$ and $v$ cannot have the same color. This completes the proof.
\end{proof}

It is straightforward that $\bar{\Gamma}$ is isomorphic to any subgraph of $\Gamma$ whose vertices form a complete representative system of the equivalence classes. From now on, we fix such a representative system and look at $\bar{\Gamma}$ as the corresponding subgraph of $\Gamma$. Doing so, $N(1)$, $N(2)$, and $N(3)$ can be restated as follows:
{\em{
\begin{itemize}
\item[N(1)] If $\bar{u}$ and $\bar{v}$ are two independent vertices of $\bar{\Gamma}$ then both $N(\bar{u})\cap N(N(\bar{v}))$ and $N(\bar{v}\cap N(N(\bar{u}))$ are the empty set.
\item[N(2)] {\mbox{$ N(N(N(\bar{u}))) \subseteq N(\bar{u})$ for each vertex $\bar{u}$ of $\bar{\Gamma}$.}}
\item[N(3)] Let $\bar{u}$ and $\bar{v}$ be two vertices of $\bar{\Gamma}$ with a common out-neighbour. If $\bar{u}\not\in N(N(\bar{v}))$ and $\bar{v}\not\in N(N(\bar{u}))$, then they have the same in-neighbours and either all out-neighbours of  $\bar{u}$ are also out-neighbors of $\bar{v}$ or all out-neighbours of  $\bar{v}$ are also out-neighbors of $\bar{u}$.
\end{itemize}
}}


\section{Paths and Circuits in a $2$-cBMG}
\label{secpc}
From now on $\Gamma=\Gamma(V,E)$ stands for a bipartite (i.e. $2$-colored) digraph without loops, parallel edges, and vertices which have no out-neighbours. Our goal is to prove some results on the paths in a 2-cBMG; see Corollary and Proposition \ref{prop11}, \ref{prop4}, and \ref{dec1}.

First, we investigate the constraints imposed on the structure of $\Gamma$ when only $N(2)$ is assumed.

\subsection{Graphs with $N(2)$}
We assume that $N(2)$ holds for $\bar{\Gamma}$. Therefore for each $u\in V$
\begin{equation}
\label{eq1} N(N(N(u))) \subseteq N(u).
\end{equation}
It may be observed that (\ref{eq1}) means that $\bar{\Gamma}$ is bitransitive. 
\begin{lemma}
\label{lem1}  Let $(\alpha_1,\alpha_2,\beta_2)$ be a triple of vertices in $\Gamma(V)$, where $\alpha_1$ and $\alpha_2$ have the same color while both $\alpha_1\beta_2$ and $\beta_2\alpha_2$ are edges in $\Gamma(E)$.
Then
\begin{itemize}
\item[(i)] $N(\alpha_2)\subseteq N(\alpha_1)$,
\item[(ii)] $N^-(\alpha_1)\subseteq N^-(\alpha_2)$.
\end{itemize}
\end{lemma}
\begin{proof} The case $\alpha_1=\alpha_2$ is trivial, therefore $\alpha_1\neq \alpha_2$ is assumed. $\alpha_1\beta_2\in \Gamma(E)$ means $\beta_2\in N(\alpha_1)$, then $N(\beta_2)\subseteq N(N(\alpha_1))$. $\beta_2\alpha_2\in \Gamma(E)$ means
 $\alpha_2\in N(\beta_2)$, this gives
$N(\alpha_2) \subseteq N(N(\beta_2))\subseteq N(N(N(\alpha_1))).$ On the other hand $N(N(N(\alpha_1))\subseteq N(\alpha_1)$ by (\ref{eq1}). Therefore $N(\alpha_2)\subseteq N(\alpha_1)$ which is claim (i).

To show claim (ii), take any vertex $\beta_1$ from $N^-(\alpha_1).$ Then $\beta_1\alpha_1\in \Gamma(E).$ This together with $\alpha_1\beta_2\in \Gamma(E)$ shows that the triple $(\beta_1,\beta_2,\alpha_1)$ satisfies the hypothesis of Lemma \ref{lem1}. From (i) applied to $(\beta_1,\beta_2,\alpha_1)$ we have $N(\beta_2)\subseteq N(\beta_1)$. Since $\alpha_2\in N(\beta_2)$, this yields $\alpha_2\in N(\beta_1)$, that is, $\beta_1\in N^-(\alpha_2)$ which proves claim (ii).
\end{proof}

\begin{lemma}
\label{lem2} Let $\alpha_1$ and $\alpha_2$ be distinct vertices of $\Gamma$ with the same color. If there exist $\beta_1,\beta_2$ (not necessarily distinct) vertices in $\Gamma$ such that
\begin{equation}
\label{eq2}
\alpha_1\beta_2\in \Gamma(E),\,\alpha_2\beta_1\in \Gamma(E),\,\beta_2\alpha_2\in \Gamma(E),\,\beta_1\alpha_1\in \Gamma(E),
\end{equation}
then $\alpha_1$ and $\alpha_2$ are equivalent.
\end{lemma}
\begin{proof}
Obviously, the color of $\alpha_1$ is different from that of $\beta_1$ and $\beta_2$. Thus $\beta_1$ and $\beta_2$ have the same color.

From the first and the third inclusion in (\ref{eq2}), Lemma \ref{lem1} holds for $(\alpha_1,\alpha_2,\beta_2)$. Therefore, we have $N(\alpha_2)\subseteq N(\alpha_1)$ and $N^-(\alpha_1)\subseteq N^-(\alpha_2)$. Using the second and the forth
inclusions, Lemma \ref{lem1} holds for ($\alpha_2$, $\alpha_1$,$\beta_1)$. Therefore, $N(\alpha_1)\subseteq N(\alpha_2)$ and $N^-(\alpha_2)\subseteq N^-(\alpha_1)$. These four inclusions together yield $N(\alpha_1)=N(\alpha_2)$ and $N^-(\alpha_1)= N^-(\alpha_2)$ which proves Lemma \ref{lem2}.
\end{proof}
\begin{lemma}
\label{lem4} Assume that $\Gamma$ contains no two equivalent vertices. Then for any $u\in \Gamma$ there exists at most one $v\in \Gamma$ such that both $uv$ and $vu$ are edges in $\Gamma(E)$.
\end{lemma}
\begin{proof}
By way of a contradiction, there are $u,v_1$ and $v_2$ vertices of $\Gamma$ such that $uv_1, v_1u, uv_2, v_2u \in E$. Note that $v_1\neq v_2$ as $\Gamma$ does not have parallel edges and that $u\neq v_1$ and $u\neq v_2$ as $\Gamma$ does not have loops. Let $w\in N(v_1)$. Then $w\in N(N(u)) \subseteq N(N(N(v_2)))$. By (\ref{eq1}), $w\in N(v_2)$. Similarly, $w\in N(v_2)$ yields $w\in N(v_1)$. Let $w\in N^{-}(v_1)$, that is $v_1\in N(w)$. Then $v_2\in N(N(N(w)))\subseteq N(w)$ whence $w\in N^{-}(v_2)$. Therefore $v_1$ and $v_2$ have the same in- and out-neighbours, contradicting our hypothesis.
\end{proof}
Lemmas \ref{lem2}  and \ref{lem4} have the following consequence.
\begin{lemma}
\label{lem3} If $\Gamma$ contains no two equivalent vertices then $\Gamma$ has no directed circuit of length four.
\end{lemma}
\begin{proof}
By way of a contradiction, let $\alpha_1\beta_1\alpha_2\beta_2\alpha_1$ be a length $4$ directed circuit of $\Gamma$. From Lemma \ref{lem4}, $\alpha_1 \neq \alpha_2$ and $\beta_1 \neq \beta_2$. Then $\alpha_1$ and $\alpha_2$, as well as $\beta_1$ and $\beta_2$, have the same color while $\alpha_1$ and $\beta_1$ have different colors. Furthermore, $\alpha_1\beta_1\in \Gamma(E)$, $\beta_1\alpha_2\in \Gamma(E)$, $\alpha_2\beta_2\in \Gamma(E)$, and $\beta_2\alpha_1\in \Gamma(E)$. Then (\ref{eq2}) holds whenever we switch $\beta_1$ and $\beta_2$, then the claim follows from Lemma \ref{lem2}.
\end{proof}

We show that Lemma \ref{lem3} is a particular case of a much stronger result.
\begin{proposition}
\label{prop11} If $\Gamma$ contains no two equivalent vertices then no directed circuits of $\Gamma$ has length grater than $2$.
\end{proposition}
\begin{proof} Since $\Gamma$ is a bipartite graph, if a directed circuit exists, it has even length. Therefore, by way of a contradiction, let $\alpha_1\beta_1\alpha_2\beta_2\cdots \alpha_i\beta_i\cdots \alpha_n\beta_n\alpha_1$ be a directed circuit of $\Gamma$ with length $2n\ge 4$. Then the vertices $\alpha_i$ have the same color, as well as the vertices $\beta_i$, where the color of $\alpha_i$ and $\beta_i$ are different. Take two consecutive vertices with the same color, say
$\alpha_i$ and $\alpha_{i+1}$. Then Lemma \ref{lem1} applies to the triple $(\alpha_i,\alpha_{i+1},\beta_i)$ showing that $N(\alpha_{i+1})\subseteq N(\alpha_i)$ and $N^-(\alpha_i)\subseteq N^-(\alpha_{i+1})$. Since this holds true for any $i$, we have
\begin{equation}
\label{eq5}
N(\alpha_n)\subseteq N(\alpha_{n-1})\subseteq \cdots \subseteq N(\alpha_2)\subseteq N(\alpha_1).
\end{equation}
Since Lemma \ref{lem1} also applies to the triple $(\alpha_n,\alpha_1,\beta_n)$ we also have $N(\alpha_1)\subseteq N(\alpha_n)$. This together with (\ref{eq5}) yields
$N(\alpha_1)=N(\alpha_2)=\cdots= N(\alpha_n)$.
Similarly, we have
\begin{equation}
\label{eq6}
N^-(\alpha_1)\subseteq N^-(\alpha_2)\subseteq \cdots \subseteq N^-(\alpha_{n-1})\subseteq N^-(\alpha_n).
\end{equation}
Applying Lemma \ref{lem1} to the triple $(\alpha_n,\alpha_1,\beta_n)$ gives $N^-(\alpha_n)\subseteq N(\alpha_1)$ which together with (\ref{eq6}) yields
$N^-(\alpha_1)=N^-(\alpha_2)=\cdots= N^-(\alpha_n)$.
Therefore, the vertices $\alpha_i$ are all equivalent, contradicting our hypothesis.
\end{proof}
The following lemma resembles \cite[Lemma 8]{geiss} for $\Gamma$.
\begin{lemma}
\label{lem24} For any two vertices $u$ and $v$ of $\Gamma,$ if $N(u)\cap N(v)=\emptyset,$ then $N(N(u))\cap N(N(v))=\emptyset.$
\end{lemma}
\begin{proof}
By way of a contradiction, there exists $t\in N(N(u))\cap N(N(v))$. Therefore, $N(t)\subseteq N(N(N(u)))\subseteq N(u)$ by (\ref{eq1}). Applying (\ref{eq1}) on $v$ and $t,$ we also have $N(t) \subseteq N(v)$. Thus, $N(u) \cap N(v)$ is not the empty if $N(t)\neq \emptyset$. Since we assumed that all vertices in $\Gamma$ have at least one out-neighbor, $N(u) \cap N(v) \neq \emptyset,$ contradicting the hypothesis.
\end{proof}
 \begin{lemma}
\label{lem21}
 For any three vertices $u,v$ and $w$ of $\Gamma$ with the same color, if $u$ and $v$ have no common out-neighbors, then either $u$ or $v$ is not strongly connected to $w$.
\end{lemma}
\begin{proof} By way of a contradiction, we have two directed walks $u=u_0\rightarrow u_1 \rightarrow u_2 \rightarrow \cdots\rightarrow u_m \rightarrow u_{m+1}=w $ and $u=v_0\rightarrow v_1 \rightarrow v_2 \rightarrow \cdots\rightarrow v_n \rightarrow v_{n+1}=w$. In particular, $u_m$ is reachable from $u$. From (\ref{eq1}), $u_m\in N(u)\cup N(N(u))$. Therefore,
either $u_m\in N(u)$ or $u_m\in N(N(u))$.  Since $u$ and $w$ have the same color different from the one of $u_m,$ the case $u_m\in N(N(u))$ cannot actually occur. Thus $u_m\in N(u)$. As $u_m w$ is an edge of $\Gamma$, it turns out that $w\in N(N(u))$. The same argument applies to $v$, therefore $w\in N(N(u))\cap N(N(v))$. But this contradicts Lemma \ref{lem24} as $u$ and $v$ have no common out-neighbors.
\end{proof}
\begin{lemma}
\label{lem22} For any three vertices $u,v$ and $w$ of $\Gamma,$ if $u$ and $v$ have the same color but have no common out-neighbors, then either $u$ or $v$ is not strongly connected to $w$.
\end{lemma}
\begin{proof} By Lemma \ref{lem21}, we may assume that $w$ does not have the same color of $u$ and $v$. Using the same argument in the proof of Lemma \ref{lem21}, $u_m\in N(u)\cup N(N(u))$. Here $u,w$ and $u_m,w$ are pairs of vertices of different colors, hence $u_m$ and $u$ have the same color and the case $u_m\in N(N(u))$ must occur. Therefore (\ref{eq1}) gives $w\in N(N(N((u)))\subseteq N(u)$. The same argument applies to $v$, then $w\in N(v)\cap N(u)$. Thus $w$ is an out-neighbor of both $u$ and $v$, contradicting one the hypotheses.
\end{proof}

\begin{lemma}
\label{lem33} Let $u$ and $v$ be any two vertices of $\Gamma$. If $u$ is out-dominated by $v,$ and $u$ is strongly connected to a vertex $w$ of $\Gamma,$ then $v$ is also strongly connected to $w$.
\end{lemma}
\begin{proof} Let $u=u_0\rightarrow u_1 \rightarrow u_2 \rightarrow \cdots\rightarrow u_m \rightarrow u_{m+1}=w$ be a directed path. Then $u_1\in N(u)$. Furthermore, $u_1\in N(v)$ since $u$ is dominated by $v$. Therefore $v=v_0\rightarrow u_1 \rightarrow u_2 \rightarrow \cdots\rightarrow u_m \rightarrow u_{m+1}=w$ is a directed path form $v$ to $w$.
\end{proof}

\subsection{Lengths of circles and paths in $2$-colored best match graphs}
\label{secdis}
Now we discuss an approach to the study of $\Gamma$ which has the advantage of connecting $2$-colored best matched graphs to more ``traditional'' graph families.

By how $\bar{\Gamma}$ arises from $\Gamma$, $\bar{\Gamma}$ is a bipartite digraph without loops, parallel edges, and vertices which have no out-neighbours. Furthermore,  $\bar{\Gamma}$ contains no two equivalent vertices but still satisfies $N(2)$. Then Proposition \ref{prop11} holds for $\bar{\Gamma}$ showing that the only directed circuits of $\bar{\Gamma}$ are those of length $2$. Therefore Proposition \ref{prop11} has the following consequence.
\begin{corollary}
\label{cor1}
The only directed cycles of $\bar{\Gamma}$ have length $2$ and are induced by symmetric edges.
\end{corollary}

Now we define a new graph $\tilde{\Gamma}$ keeping the same vertex set but eliminating all symmetric edges from $\bar{\Gamma}$. In other words, if for some $u,v\in \bar{\Gamma}(\bar{V})$ both $uv$ and $vu$ are edges in $\bar{\Gamma}$, then we keep only one of them. From Corollary \ref{cor1}, the oriented bipartite graph $\tilde{\Gamma}$ with no isolated vertex is acyclic.
However we may have some vertex without out-neighbour.
Result \ref{dastheorem3.4} shows that acyclic oriented digraphs have a well-determined structure. Such digraphs also pay a role in genomics data processing; see \cite{appl}. 

The results reported in Section \ref{secpc} show that 
under natural hypotheses 
$\tilde{\Gamma}$ has longer directed paths. For instance, Result \ref{zanghA} yields that this occurs when $k+h$ is big enough. In \cite[Fig. 7 A]{geiss}, whenever we derive $\tilde{\Gamma}$ from $\bar{\Gamma}$, we end up with two vertices one with no in-neighbours and the other with no out-neighbour, that is $k,h=0$. The above result ensures the existence of a directed path of length $3$ in $\tilde{\Gamma}$ and hence in $\bar{\Gamma}$. Actually $\bar{\Gamma}$ contains the directed path $\alpha_1\rightarrow \beta_1 \rightarrow \alpha_5  \rightarrow \beta_3 \rightarrow \alpha_3$ of length $4$.

Several other results about bipartite oriented graphs with longer paths are available in the literature, here we limit ourselves to cite some of them.

\subsection{Graphs with $N(1)$ and $N(2)$}
\label{n1n2}
In this section we assume Property $N(1)$ and Property $N(2)$. Then (\ref{eq1}) holds for $\Gamma$, and for any two distinct vertices $u,v$ of $\Gamma$
\begin{equation}
    \label{n1}
    \mbox{$u\not\in N(v)$ or $v\not\in N(u)$ implies $N(N(v))\cap N(u)=N(v)\cap N(N(u))=\emptyset.$}
\end{equation}
\begin{lemma}
\label{lem23} For any three vertices $u,v$ and $w$  of $\Gamma$, if $u$ and $v$ have different colors but they are independent, then either $u$ or $v$ is not strongly connected to $w$.
\end{lemma}
\begin{proof}
By way of a contradiction, both $u$ and $v$ are assumed to be strongly connected to $w$. W.l.g. we may suppose that $w$ and $v$ have the same color different from the one of $u$. The argument in the proof of Lemma \ref{lem22} applied to $u,w$ shows that $w\in N(u),$ while the argument in the proof of Lemma \ref{lem21} applied to $v,w$ gives $w\in N(N(v))$. Therefore $N(u)\cap N(N(v))\neq\emptyset$. Then $u$ and $v$ are not independent by (\ref{n1}), contradicting one of the hypotheses.
\end{proof}
Lemma \ref{lem23} has the following consequence.
\begin{proposition}
\label{prop4} Let $u$ and $v$ be two independent vertices of $\Gamma$ with no common out-neighbors. Then, for any vertex $w$, either $u$ or $v$ is not strongly connected to $w$.
\end{proposition}

\subsection{Graphs with $N(2)$ and $N(3)$ }
\label{n2n3}
Finally we assume that  $N(2)$ and $N(3)$ hold for $\bar\Gamma$. Then (\ref{eq1}) holds and
 for any two different vertices $u$ and $v$ with a common out-neighbour in $\Gamma$ such that $u\not\in N(N(v))$ and $v\not\in N(N({u}))$ we have
\begin{equation}
    \label{n3}
    \mbox{ $N^{-}(u)=N^{-}(v)$ and  $u$ is dominated by $v$ or $v$ is dominated by $u$.}
\end{equation}

\begin{proposition}
\label{dec1}
Let $u$ and $v$ be non-equivalent vertices in $\Gamma$ with a common out-neighbour. If there is no a directed path of length $2$ from $u$ to $v$ or vice-versa, then at least one of them is not the endpoint of a symmetric edge.
\end{proposition}
\begin{proof}
By our hypotheses, $u$ and $v$ satisfy (\ref{n3}). W.l.g. $u$ is dominated by $v$. Assume that there exists a vertex $w$ of $\Gamma$ such that $uw,wu\in E$. Then $wv\in E$ by $N^{-}(u)=N^{-}(v)$. Furthermore $vw\in E$ since $u$ is dominated by $v$. Now the claim follows from Lemma \ref{lem4} applied to $w$.
\end{proof}
In \cite[Fig.7 A]{geiss}, $N(\alpha_5)=\{\beta_2,\beta_3,\beta_4\}$ and $N(\alpha_6)=\{\beta_3,\beta_4\}$, $N(\alpha_5)=\{\alpha_2,\alpha_3,\alpha_4\}$, and $N(\alpha_6)=\{\alpha_3,\alpha_4\}$. Also, $\alpha_6$ is dominated by $\alpha_5$. This shows the hypotheses of Proposition \ref{dec1} are satisfied by $u=\alpha_6$ and $v=\alpha_5$. Accordingly, either $\alpha_5$ or $\alpha_6$ is not the endpoint of symmetric edge. Here actually both have the property.

\section{Paths and Hierarchy}

As above let $\Gamma$ denote a connected bipartite graph without loops, parallel edges, and vertices which have no out-neighbour. From \cite[Lemma 9]{geiss} if $\bar{\Gamma}$ satisfies properties $N(1),$ $N(2)$ and $N(3)$, then it has a hierarchy property. This means that the set $R(\bar{\alpha})=N(\bar\alpha)\cup N(N(\bar\alpha))$ of all vertices of $\bar{\Gamma}$ reachable from $\alpha$ has the following properties. For any two vertices $\bar{\alpha}$ and $\bar{\beta}$, either  $R(\bar{\alpha})\subseteq R(\bar{\beta})$, or $R(\bar{\beta})\subseteq R(\bar{\alpha})$, or $R(\bar{\alpha})\cap R(\bar{\beta})=\emptyset$. Comparing this property to the results in Section \ref{secpc} suggests that some of those results, for instance, Lemma \ref{lem22}, can be reproved by using the hierarchy property.

\end{document}